\documentclass
[12pt]{amsart}

\usepackage{amssymb}

\usepackage[all]{xy}

\newcommand{\A}{{\mathbb{A}}}
\newcommand{\R}{{\mathbb{R}}}

\newcommand{\C}{{\mathbb{C}}}

\newcommand{\N}{{\mathbb{N}}}

\newcommand{\GL}{\operatorname{GL}}

\newcommand{\rank}{{\rm rank}}

\theoremstyle{plain}
\newtheorem{thm}{Theorem}[section]
\newtheorem{lem}[thm]{Lemma}
\newtheorem{prop}[thm]{Proposition}
\newtheorem{cor}[thm]{Corollary}

\theoremstyle{definition}

\address{Institut Camille Jordan, Universit\'e Claude Bernard - Lyon
I, B\^atiment de Math\'ematiques, 43, Bld. du 11 Novembre 1918,
69622 Villeurbanne Cedex, France {\tt tomanov@math.univ-lyon1.fr}}

\begin{document}

\evensidemargin 2cm
\oddsidemargin 2cm

\centerline{\bf {PROPERLY DISCONTINUOUS GROUP ACTIONS ON}}
\centerline{\bf {AFFINE HOMOGENEOUS SPACES}}
\vskip.7cm

\centerline{George Tomanov}
\vskip.4cm


\vskip.6cm

\section{Introduction} \label{Introduction}

Let $G$ be a real algebraic group, $H$ an algebraic subgroup of $G$,
and $\Gamma$ a closed subgroup of
$G$ acting on the
homogeneous space $G/H$ by left translations. Given $x\in G/H$, $\Gamma_x$ is the stabilizer of $x$ in $\Gamma$. Recall that the action
of $\Gamma$ is {\it properly discontinuous} (respectively, {\it free})
if
for any compact $K \subset G/H$ the set
$\{ g \in \Gamma \mid gK \cap K \neq \emptyset \}$ is finite
(respectively,  $\Gamma_x$ is trivial for all $x\in G/H$). If $\Gamma$ acts properly
discontinuously and freely on $G/H$ then the manifold of double co-sets
$\Gamma \backslash G/H$ is called  Clifford-Klein form.
The following question is natural and well-known: Which homogeneous manifolds $G/H$
admit nontrivial (respectively, compact) Clifford-Klein forms
 $\Gamma \backslash G/H$? The question has been studied when the homogeneous spaces  $G/H$ is of reductive
type, that is, when both $G$ and
$H$ are reductive groups (cf.[Be1-2] and [K1-3]). In the present
paper  we discuss some cases when $G/H$ is never of
reductive type. First of all recall the notable
Auslander conjecture  (cf. [Au]):

\vskip.5cm
{\bf Conjecture 1}. {\it Let $\Gamma$ be a subgroup
of the group $\mathrm{Aff}({\mathbb{A}}^n)$ of all affine linear transformations of the
$n$-dimensional real affine space ${\mathbb{A}}^n$. Assume that $\Gamma$ acts properly
discontinuously on ${\mathbb{A}}^n$ and the quotient $\Gamma \backslash
{\A}^n$
is compact. Then $\Gamma$ is a virtually solvable group, i.e. $\Gamma$
 contains a solvable subgroup of finite index. }
\vskip.5cm

In other terms, the Auslander conjecture says that if $\Gamma$ acts properly
discontinuously and co-compactly on ${\mathbb{A}}^n$ then the Levi subgroup of its Zariski closure in $\mathrm{Aff}({\mathbb{A}}^n)$
is trivial. (Recall that the maximal connected semisimple subgroups of $G$ are usually called Levi subgroups of $G$ and they are all conjugated.)
As proved by G.A.Margulis in \cite{[Mar1]} and \cite{[Mar2]}, the compactness of $\Gamma \backslash{\A}^n$ in the formulation of the conjecture is essential. (See Theorem 4.2(b) below.)

The continuous analog of the Auslander conjecture  is the
following result of
T.Kobayashi and R.Lipsman:

\begin{thm} \mbox{ \rm{(}} \cite{[K1]}, \cite{[L]} \mbox{\rm{)}}
\label{KL}
Suppose that $H$ contains a Levi subgroup of $G$, $\Gamma$ is a
connected algebraic subgroup of $G$ and
$\Gamma_x$ is compact for all $x \in G/H$. Then $\Gamma$
is a compact extension of a unipotent group.
\end{thm}

In the light of the above discussion the following generalization of Auslander's conjecture is
natural:

\vskip.5cm
{\bf Conjecture 2.} {\it Suppose that $H$ contains a maximal reductive subgroup of $G$,
$\Gamma$ acts properly discontinuously on $G/H$ and $\Gamma \backslash G/H$ is compact. Then $\Gamma$
is virtually solvable.}
\vskip.5cm

It is easy to see that that Conjecture 2 implies Conjecture 1 (cf. Remark 1 in 2.1). Also note that  $G/H$ is isomorphic (as a real
algebraic variety) to ${\A}^n$ and $G$ acts on ${\A}^n$ by regular (polynomial) automorphisms of degree $\geq 1$. Conjecture 1 is exactly the case when
this action is linear.
Some known results about Auslander's conjecture could be extended to Conjecture 2.
In \S2 and \S3 of the present paper we prove the following
\begin{thm}
\label{Conj2} Conjecture 2 is true if the Levi subgroup of
$G$ is a product of simple real algebraic groups of ranks $\leq 1$.
\end{thm}

The arguments used in the proof of Theorem \ref{Conj2} generalize
our arguments in \cite{[To1]} where the analogous result is proved for Auslander's conjecture\footnote{Our result \cite{[To1]} can be
also found in \cite{[So1]}.}.
Independently, K.Dekimpe and N.Petrosyan proved in a recent paper a similar to Theorem 1.2 result \cite[Theorem A]{[D-P]} and formulated relevant to  Conjecture 2 questions \cite[Questions 1 and 2]{[D-P]}. 
At present, we can prove Conjecture 2 for $\dim G/H \leq 4$. (The proof will appear elsewhere.)

As to the Auslander conjecture, its proof is easy for $n = 2$
and due to D.Fried and W.Goldman \cite{[F-G]} for $n = 3$. The proof of the conjecture
for $n \leq 5$ was announced in \cite{[To2]} with detailed sketch of the proof for $n = 4$.
After the present paper was finished we became
aware of the preprint \cite{[A-M-S5]} where the conjecture is proved for $n \leq 6$. In \S4 we
give a full proof of Auslander's conjecture for $n \leq 5$ (Theorem \ref{thm small dimensions}) along the lines in \cite{[To2]}.
Our proof for $n \leq 5$ is different and simpler than the proof in \cite{[A-M-S5]}, in particular, it uses less input. All one needs is the result of
Margulis preprint \cite{[Mar3]}, published as part of \cite{[A-M-S3]} (see Theorem \ref{AMS3}(a) below), and \cite{[To1]}.

\subsection{Notation and terminology.}\label{notation} By an algebraic group (resp., algebraic
variety) we will mean a {\it real} linear algebraic group (resp., {\it real}
algebraic variety), that is, the set of all ${\R}$-rational points of a linear
algebraic group (resp., algebraic variety)
defined over ${\R}$. On every algebraic variety we have Hausdorff topology (induced by the
topology on $\R$) and also Zariski topology. In order to distinguish the two topologies the topological notions
connected with the Zariski topology will be usually used with the prefix "Zariski". (We say: Zariski closed, Zariski closure, Zariski connected, etc.) If $M$ is a subset of an algebraic variety $X$
then $\tilde M$ denotes the Zariski closure of $M$ in $X$ and $\overline{M}$ denotes the closure of $M$ in $X$ for
the usual (Hausdorff) topology.  
We denote by $G^\circ$ the connected component of $G$ for the Hausdorff topology and by $\mathrm{R}(G)$ (respectively, $\mathrm{R}_u(G)$)
the radical (respectively, unipotent radical) of $G$. If $G$ acts on a set $X$ and $x \in X$, $G_x$ is
the stabilizer of $x$ in $G$. Given $g \in G$, $g = g_s g_u$ is the Jordan decomposition of $g$
where $g_s$ (resp. $g_u$) is the semi-simple (resp. unipotent) part of $g$. We let $<g>$ be the subgroup generated by $g$. Also, we will
denote by $\mathrm{Lie}(G)$ the Lie algebra of $G$. By \textit{rank of} $G$ we mean
the common dimension of the maximal
$\R$-diagonalizable tori of $G$.
Also,  ${\mathcal D}^i G$ is the $i$-th derived subgroup of $G$, that is, $\mathcal{D}^0 G = G$ and ${\mathcal D}^{i+1} G = [{\mathcal D}^i G, G]$ for all $i \geq 0$.
\subsection{Basic affine geometry.}\label{affine geometry}
A real affine space $\A^n$ is obtained from a real $n$-dimensional vector space $V$, called the \textit{direction} of $\A^n$, by "forgetting"
the origin. Most often $V = \R^n$.
If $x$ and $y \in \A^n$ we denote by $\overrightarrow{xy}$ the unique vector in $V$ such that $y = x + \overrightarrow{xy}$. An affine automorphism $\gamma \in \mathrm{Aff}(\A^n)$ determines a $\lambda(\gamma) \in \GL(V)$, called  \textit{the linear part} of $\gamma$, such that for any pair $x, y \in \A^n$ we have $\overrightarrow{\gamma(x)\gamma(y)} = \lambda(\gamma)(\overrightarrow{xy})$. The map $\lambda: \mathrm{Aff}(\A^n) \rightarrow \mathrm{GL}(V), \gamma \mapsto \lambda(\gamma)$, is a surjective group homomorphism. Fix an \textit{origin} $p \in \A^n$.
If $\overrightarrow{v} \in V$ then $\gamma(p + \overrightarrow{v}) = \gamma(p) + \lambda(\gamma)(\overrightarrow{v})$. So, every $\gamma \in \mathrm{Aff}(\A^n)$ can be decomposed as the linear "vector" transformation $\A^n \rightarrow \A^n, p + \overrightarrow{v} \mapsto p + \lambda(\gamma)(\overrightarrow{v})$, followed by the translation by the
vector $\overrightarrow{p\gamma(p)}$. Consider the semidirect product of algebraic groups $V \rtimes \mathrm{GL}(V)$ where the action of $\mathrm{GL}(V)$ on $V$ is the natural one. The group $\mathrm{Aff}(\A^n)$ is identified with $V \rtimes \mathrm{GL}(V)$ via the group isomorphism $\mathrm{Aff}(\A^n) \rightarrow V \rtimes \mathrm{GL}(V), \gamma \mapsto (\overrightarrow{p\gamma(p)}, \lambda(\gamma))$.
The structure of algebraic group on $\mathrm{Aff}(\A^n)$ obtained in this way does not depend on the choice of $p$.
All stabilizers ${\mathrm{Aff}(\A^n)}_x, x \in \A^n,$
are maximal reductive
subgroups of $\mathrm{Aff}(\A^n)$ isomorphic to $\mathrm{GL}(V)$ and pairwise vector translation conjugate.
Hence every reductive subgroup of $\mathrm{Aff}(\A^n)$ (and, therefore, every semi-simple element in $\mathrm{Aff}(\A^n)$) admits a fixed point. Further on, we will tacitely use this observation.

The group $\R^n\rtimes\mathrm{GL}_n(\R)$ is identified with its image in $\GL_{n+1}(\R)$ under the imbedding $(\overrightarrow{v}, m) \mapsto \left( \begin{array}{cc}  m&\overrightarrow{v}  \\  0&1 \\  \end{array}
 \right)$, where the elements from $\R^n$ are vector columns.
 If $\overrightarrow{v_1}, \cdots, \overrightarrow{v_n}$ is a basis of $V$, then $\mathcal{F} = \{p; \overrightarrow{v}_1, \cdots, \overrightarrow{v}_n\}$ is a \textit{frame} of $\A^n$. We have an isomorphism $\mathrm{Aff}(\A^n) \rightarrow \R^n\rtimes\mathrm{GL}_n(\R), \gamma  \mapsto \left( \begin{array}{cc}  l(\gamma)&\overrightarrow{v}(\gamma)  \\  0&1 \\  \end{array}
 \right)$, where $l(\gamma)$ is the matrix of $\lambda(\gamma)$ in the basis $\overrightarrow{v_1}, \cdots, \overrightarrow{v_n}$ and $\overrightarrow{v}(\gamma)$ is the vector-column of the coordinates of $\overrightarrow{p\gamma(p)}$ in this basis, called
  \textit{the matrix representation of} $\mathrm{Aff}(\A^n)$ \textit{in the frame} $\mathcal{F}$.  Note that $l(\gamma)$ is the
  same in any translated frame $\mathcal{F} + \overrightarrow{v} = \{p+ \overrightarrow{v}; \overrightarrow{v}_1, \cdots, \overrightarrow{v}_n\},
  \overrightarrow{v} \in V$.

\section{Rational actions of $\Gamma$ on $\A^n$} \label{Proof thm}

\subsection{}
Let $G$ and $\Gamma$ be as in the formulation of Conjecture 2. Replacing $\Gamma$ by a subgroup of finite index, we suppose
  from now on that $G$ is Zariski connected. Our main goal is to reduce the proof of
Conjecture 2 to the case when $\Gamma \cap \mathrm{R}(G) = \{ e \} $.
\vskip.5cm

\begin{prop}
\label{radical}
 Assume that $G$ is acting rationally on $\A^n$, the restriction of this action to $\Gamma$
is properly discontinuous, the quotient $\Gamma \backslash \A^n$ is compact, and there exists $x_o \in \A^n$ such that $S x_o = x_o$ for a maximal reductive subgroup $S$ of $G$.
Then $\mathrm{R}_u(G)$
acts transitively on $\A^n$.
\end{prop}

{\bf Proof.} Note that
$Gx_o = \mathrm{R}_u(G)x_o$ and $\mathrm{R}_u(G)x_o$ is closed and isomorphic as a real algebraic variety to an affine space $\A^k$ (see \cite{[Bi]} and \cite{[Ro]}).
The group
 $\Gamma$ acts properly discontinuously and with compact quotient
on $\mathrm{R}_u(G)x_o$.  In view of [Se],
$$
\mbox{vcd}(\Gamma) = \dim \mathrm{R}_u(G)x_o = \dim {\A}^n,
$$
where vcd$(\Gamma)$ denotes the virtual cohomological dimension of
$\Gamma$.  So, $\mathrm{R}_u(G)x_o = \A^n$, completing the proof. \qed

\vskip.5cm
{\it Remarks:} 1. Let $\Gamma$ be as in the formulation of the
Auslander conjecture, $G$ be the Zariski closure of $\Gamma$ in
$\mbox{Aff}({\A}^n)$ and $S$ be a maximal reductive subgroup of
$G$. Let $x_o \in {\A}^n$ be fixed by $S$.
Put $H = G_{x_\circ}$. In view of
Proposition \ref{radical}, $G$ acts transitively
on ${\A}^n$. So,  ${\A}^n$ can be identified with $G/H$ and, therefore,
 Conjecture 2 implies Conjecture 1.

 2. The argument used in the proof of Proposition \ref{radical} is identical to that used in the proof of \cite[Lemma 1.1]{[To1]}.
As indicated to the author by the referee of \cite{[To1]}, in a different way, the result
 was first proved by W.Goldman and M.W.Hirsch \cite[Theorem 2.6]{[G-H]}.  \cite[Lemma 2.5]{[D-P]} corresponds to
 \cite[ Lemma 1.1]{[To1]}.

\medskip

 The following  assertion is implicitly contained in the proof of
\cite[Proposition 1.4]{[To1]}.
\vskip.5cm

\begin{lem}
\label{solv}
 Let $\Delta \subset {\mathrm{GL}_n}({\R})$ be a discrete
solvable subgroup. Then there exists a connected (for the Hausdorff topology on ${\mathrm{GL}_n}({\R})$) solvable subgroup $R
\subset {\mathrm{GL}_n}({\R})$ such that $R \cap \Delta$ is a normal subgroup of finite index in
$\Delta$ and $R/R \cap \Delta$ is compact.
\end{lem}
{\bf Proof.} The group $\Delta \cap \tilde \Delta^\circ$ is a normal subgroup of finite index in
$\Delta$. Replacing $\Delta$ by $\Delta \cap \tilde \Delta^\circ$ it is enough to prove the existence of a connected
subgroup $R$ such that $\Delta \subset R \subset \tilde \Delta^\circ$ and $R/\Delta$ is compact.
 Note that ${\mathcal D} \Delta$ is a Zariski dense
discrete subgroup in the connected unipotent group ${\mathcal D} \tilde \Delta = {\mathcal D} \tilde \Delta^\circ$.
Therefore $\Delta \cap {\mathcal D} \tilde \Delta$ is a co-compact lattice in ${\mathcal D} \tilde \Delta$.
This implies that $\Delta \cdot {\mathcal D} \tilde \Delta/{\mathcal D} \tilde \Delta$ is a discrete subgroup of $\tilde \Delta^\circ/{\mathcal D} \tilde \Delta$
(cf. \cite[Theorem 1.13]{[Rag]}). Hence, it is enough to prove the lemma when $\tilde \Delta$ is abelian and $\Delta \subset \tilde \Delta^\circ$.
 The Lie group $\tilde \Delta^\circ$ is isomorphic to $K \times \R^m$ where $K$ is a compact torus. Let $\pi: \tilde \Delta^\circ \rightarrow \R^m$ be the natural projection and $R'$ be the linear span of $\pi(\Delta)$. Then $\Delta$ is co-compact in $R = \pi^{-1}(R')$. \qed

\subsection{}
Let $G, H$ and $\Gamma$ be as in the formulation of
Conjecture 2. Using Proposition \ref{radical} we see that the unipotent radical of the Zariski closure of $\Gamma$ in $G$ is acting
transitively on $G/H$. So, replacing $G$ by the Zariski closure of $\Gamma$ we may (as we will) assume
 that $\Gamma$ is Zariski dense in $G$. 
 Put $\Delta = \Gamma \cap \mathrm{R}(G)$.
Then $\tilde \Delta$ is a normal subgroup
of $G$. Denote by $G_1$ the Zariski closure of $G/\tilde \Delta$, by $H_1$ the Zariski closure of $H\tilde
\Delta /\tilde \Delta$ in $G_1$, and by $\Gamma_1$ the natural
imbedding of $\Gamma/\Delta$ into $G_1$. Clearly, $\Gamma_1 \cap
\mathrm{R}(G_1) = \{e\}$.

The next proposition, which is the central one, allows
to "eliminate" the solvable radical when dealing with the Auslander conjecture or with some of its generalizations. Actually, it coincides with \cite[Proposition 1.4(a)]{[To1]}.
 For reader's convenience we provide a somewhat more detailed than in \cite{[To1]} proof of the proposition.

\vskip.5cm

\begin{prop}
\label{reduction}
With the above notation and assumptions,
$\Gamma_1$ acts properly discontinuously on $G_1/H_1$ and
$\Gamma_1\backslash G_1/H_1$ is compact.
\end{prop}

{\bf Proof.} Since $\tilde \Delta$ is a normal subgroup in $G$,  the
action of $G$ on $G/H$ permutes the $\tilde \Delta$-orbits on
$G/H$. So, we can identify the space of $\tilde
\Delta$-orbits on $G/H$ with $G/H'$ where $H' = \tilde \Delta
H$. Let $\tilde \Delta_u$ be the unipotent radical of
$\tilde \Delta$ and $T$ be a maximal reductive subgroup of $\tilde
\Delta$. Then $\tilde \Delta$ is equal to the semidirect product $\tilde \Delta_u \rtimes T$. Remark that $T$ is conjugated to a subgroup of $H$
and $\tilde \Delta_u$ is normal in $G$.
Hence $H' = \tilde \Delta_u H$ and $H'$ is an algebraic subgroup of $G$.

Let $\phi:\tilde \Delta \rightarrow \tilde \Delta_u$
be the natural  projection and $R$ be a connected subgroup of $\tilde \Delta$
such that $\Delta \cap R$ is a normal subgroup of finite index in $\Delta$ and $R/\Delta \cap R$ is compact (see Lemma \ref{solv}).
We will prove that $\phi(R) = \tilde \Delta_u$. Denote by $\tilde{\Delta}^\bullet$ the Zariski connected component of $\tilde \Delta$.
Then $\tilde \Delta_u$ is the unipotent radical of $\tilde{\Delta}^\bullet$. Suppose that $\tilde{\Delta}^\bullet$ is abelian. In this case
the restriction of $\phi$ to $\tilde{\Delta}^\bullet$ is a homomorphism of algebraic groups and $\phi(R)$ is connected and , therefore, algebraic subgroup of $\tilde{\Delta}_u$. Since $\phi(R)$ is Zariski dense in $\tilde{\Delta}_u$ we get that $\phi(R) = \tilde \Delta_u$. Now, let $\tilde{\Delta}^\bullet$ be arbitrary. Since $R$ is connected, the commutator $\mathcal{D}(R)$ is unipotent and $\phi(R)$ contains  $\mathcal{D}(R)$. It is enough to prove that  $\mathcal{D}(R) = \mathcal{D}(\tilde{\Delta}^\bullet)$.
Indeed, if so, we may factorize by  $\mathcal{D}(R)$ and reduce the proof to the case when  $\tilde{\Delta}^\bullet$ is abelian. Let us prove
that $\mathcal{D}(R)$ contains $\mathcal{D}(\tilde{\Delta}^\bullet)$. (The inclusion $\mathcal{D}(R) \subset \mathcal{D}(\tilde{\Delta}^\bullet)$ is obvious.)
Since $R$ is Zariski dense in $\tilde{\Delta}^\bullet$
and $\mathcal{D}(R)$ is an algebraic subgroup of  $\tilde{\Delta}^\bullet$ we have that $\mathcal{D}(R)$ is normal in $\tilde{\Delta}^\bullet$. But $R/\mathcal{D}(R)$ is Zariski dense in  $\tilde{\Delta}^\bullet/\mathcal{D}(R)$. Therefore $\tilde{\Delta}^\bullet/\mathcal{D}(R)$ is abelian which implies that $\mathcal{D}(R)$ contains $\mathcal{D}(\tilde{\Delta}^\bullet)$, as required.

In view of the above, if $x \in G/H$ then $\tilde \Delta x = \tilde \Delta_u x = Rx$ is closed and the quotient
$\Delta\backslash \tilde \Delta x$ is compact. Since $\Delta$ acts trivially on $G/H'$, the natural action of
$\Gamma$ on $G/H'$ induces an action of $\Gamma_1$ on $G/H'$. Let us prove that
$\Gamma_1$ acts properly discontinuously on $G/H'$ and that $\Gamma_1\backslash G/H'$ is compact. Indeed, let
$\psi:G/H \rightarrow G/H'$ be the natural map, $K_o \subset G/H$ be a
compact subset and $K = \psi(K_o)$. Since $\psi$ is $\Gamma$-equivariant we have that $\Gamma_1 K =  G/H'$ if
$\Gamma K_o =  G/H$, proving that $\Gamma_1\backslash G/H'$ is compact. Let $\{{\gamma_i}' \mid
i\in I \}$ be the set of all elements in $\Gamma_1$ such that
$ {\gamma_i}' K\cap K \neq \emptyset $. For each $i$ we fix a
$\gamma_i \in \Gamma$ such that ${\gamma_i}' = \gamma_i
\Delta$.
Every fiber of $\psi$ is a $\tilde \Delta$-orbit and, by the above,
an $L$-orbit. Therefore for every $i \in I$ there
exist $a_i, b_i \in K_o$ and $l_i \in R$ such that $\gamma_ia_i = l_ib_i$.
Fix a compact $C \subset R$ such that $R = \Delta C$ and write $l_i =
\delta_ic_i$, where $ \delta_i \in \Delta$ and $ c_i \in C$. Then
$(\delta^{-1}_{i} \gamma_i)a_i = c_ib_i$. But $\Gamma$ acts properly
discontinuously on $G/H$. Therefore $\{\delta_{i}^{-1} \gamma_i \mid i
\in I \}$ is finite, which implies that $\Gamma_1$ acts properly
discontinuously on $G/H'$.

In order  to complete the proof of the proposition it remains to
notice that  $G/H'$ and $G_1/H_1$ are both canonically homeomorphic to
the affine variety $\mathrm{R}_u(G)/\widetilde{\Delta}_u \cdot \mathrm{R}_u(H)$.
\qed

We will use Proposition \ref{reduction} together with the following:

\begin{prop}
\label{reduction1} With the notation and assumptions of Conjecture 2,
additionally assume that $\Gamma \cap \mathrm{R}(G) = \{e\}$ and $\Gamma$ is Zariski dense in $G$.
Denote  by $L$ a Levi subgroup of $G$ and by $K$ a maximal compact
subgroup of $L$. Then
\begin{equation}
\label{1}
\dim(G/H) \leq \dim(L/K).
\end{equation}
\end{prop}
{\bf Proof.} Since $\Gamma$ acts properly discontinuously and with
compact quotient on the affine space $G/H$, we have
$$
\dim(G/H) = \mbox{vcd}(\Gamma).
$$
On the other hand, the projection of $\Gamma$ into $G/\mathrm{R}(G)$ is
injective and the connected component of its closure in $G/\mathrm{R}(G)$ is solvable by a result of
Auslander (cf.[Rag, Theorem 8.24]).
Therefore the projection of $\Gamma$ into $G/\mathrm{R}(G)$ is discrete. This implies that $\Gamma$ acts
properly discontinuously on the symmetric space of $L$. Therefore,
$$\mbox{vcd}(\Gamma) \leq L/K, $$
completing the proof. \qed

\medskip

The following proposition is useful.

\begin{prop}
\label{semisimple element}
With $G$ and $H$ as in the formulation of Conjecture 2, let $g \in G$. Let $U = \widetilde{<g_u>}$ where $g_u$ is the unipotent part of $g$. Then there
exists $p \in G/H$ with the following properties:
\begin{enumerate}
\item[(i)] The orbit $Up$ is closed and $g$-invariant;
\item[(ii)] $g_s$ fixes $Up$ element-wise;
\item[(iii)] $\dim Up = 1$ if $g_u \neq e$ and $gp = p$ if $g_u = e$.
\end{enumerate}
\end{prop}
{\bf Proof.} Since $H$ contains a maximal reductive subgroup of $G$ there exists a $\sigma \in G$ such that
$g_s \in \sigma H\sigma^{-1}$. Let $p = \sigma H$. It follows from $g_s g_u = g_u g_s$ that $g_s$ fixes $Up$ element-wise and that $Up$
if $g$-invariant. It is well known (and easy to prove) that  $\dim U = 1$ if $g_u \neq e$. Finally, $Up$ is closed as a unipotent orbit
on an affine algebraic variety (cf.\cite{[Bi]}).
\qed

\medskip

{\it Remark:}
If $g \in \mathrm{Aff}(\A^n)$ and
all eigenvalues of $\lambda(g)$ are different from $1$ and pairwise different then $g = g_s$ and according to (iii) (applied to $G = \mathrm{Aff}(\A^n)$ and $H = \GL_n(\R)$) there exists a $p \in \A^n$ such that $g p = p$.
This assertion is well-known and is easy to prove directly. Using it one proves easily that if $\Gamma$ acts properly discontinuously on $\A^2$ then $\Gamma$ is virtually solvable, in particlar, the Auslander conjecture holds when $n = 2$.

\medskip

Finally, let us also mention:

\begin{prop}
\label{linearization} Let $G$ and $H$ be as in the formulation of Conjecture 2 and $S$ be a maximal reductive subgroup of $G$ contained in $H$.
 Then the action of $S$ on $G/H$ by left
translations is linearizable.
\end{prop}

{\bf Proof.} Put $U=\mathrm{R}_u(G)$.  Then $U_1 = H
\cap U$ is the unipotent radical of $H$ and $G/H$ is rationally
isomorphic to $U/U_1$. Since $U$ and $U_1$ are
$\mathrm{Int}(S)$-invariant there exists $\mathrm{Ad}(S)$-invariant vector subspace
${\mathcal W} \subset
\mbox{Lie}(U)$ such that Lie$(U) = {\mathcal W} \bigoplus \mbox{Lie}(U_1)$ and
$W = $ exp${\mathcal W}$ is a
regular cross-section for $U/U_1$ (i.e. the map $ W \times U_1
\rightarrow U$, $(x,y) \rightarrow xy$, is a regular isomorphism of real
algebraic varieties), cf. [Bo-Spr, 9.13]. Since
$\exp \circ \mbox{Ad}(x) = $Int$(x) \circ $exp for any $x \in S$, we have that
the map ${\mathcal W} \rightarrow G/H$, $ w \rightarrow (\exp w)H$, is
$S$-equivariant isomorphism of algebraic varieties. Therefore, the
action of $S$ on $G/H$ is linearizable. \qed

\begin{cor}
\label{linearization+} Suppose that the action of $S$
on $G/H$ is irreducible. Then the action of $G$
on $G/H$ is linearizable, that is, there exists an isomorphism $\varphi: G/H \rightarrow \A^n$
such that if $g \in G$ and $l_g: G/H \rightarrow G/H, xH \mapsto gxH,$ then $\varphi \circ l_g \circ \varphi^{-1} \in \mathrm{Aff}(\A^n)$
for all $g$.
\end{cor}
{\bf Proof.} We use the notation $U$ and $U_1$ as in the proof of Proposition \ref{linearization}. Let $\mathcal{N}_U(U_1)$ be the normalizer
of $U_1$ in $U$. We suppose that $\dim G/H > 0$. Then $\mathcal{N}_U(U_1) \supsetneqq U_1$. Since $\mathcal{N}_U(U_1)$ is $S$-invariant
and the action of $S$ on $G/H$ is irreducible, $U_1$ is a normal subgroup of $U$ and, therefore, of $G$ too. Factorizing $G$ and $H$ by $U_1$ we reduce the proof to the case when $U_1 = \{e\}$, i.e., when $H = S$. Since $\mathcal{D}(U) \cdot S$ is a proper subgroup of $G$ and the action of $S$ on $G/H$ is irreducible, $\mathcal{D}(U)= \{e\}$. Hence $G$ is a semidirect product of $S$ and the vector group $U$ on which $S$ acts linearly, implying the corollary. \qed

\medskip

 Corollary \ref{linearization+} shows that Conjectures 1 and 2 coincide for irreducible actions of $S$ on $G/H$.

\section{Proof of Theorem \ref{Conj2}}

\subsection{Some representation theory.} Let $G$, $H$, $\Gamma$ be as in the formulation of Theorem \ref{Conj2}.
Let $L$ be a Levi subgroup of $H$. Then $L$ is an almost direct product of
simple algebraic groups $L_1, \cdots, L_r$ each of rank $\leq 1$.

Using Proposition \ref{reduction} we reduce the proof of the theorem
to the case when $\Gamma \cap \mathrm{R}(G)$ is trivial and $\Gamma$ is Zariski dense in $G$.
Moreover, by a theorem of Selberg (see \cite{[S]}) $\Gamma$ contains a torsion free subgroup of finite index.
Hence, we may (and will) suppose that $\Gamma$ is torsion free.
In view of Proposition \ref{reduction1}
the relation (\ref{1}) holds. We will denote by $V$ the tangent space of $G/H$ at the origin and
by $\rho$ the representation of $L$ on $V$ induced by the action of $L$
on $G/H$ by left translations (Proposition \ref{linearization}). Since the kernel of the action of $G$ on $G/H$ is a normal
algebraic subgroup $N$ of $G$ contained in $H$, factorizing $G$ and
$H$ by $N$ we may (and will) suppose that $G$ acts faithfully on
$G/H$. In this case the representation $\rho$ is also faithful.

The following proposition is an improved version of \cite[2.5]{[To1]}.

\begin{prop}
\label{representation} 
  With the above notation and
assumptions, $L=S_1 \times S_2 \times ... \times S_m$, where $S_i = \mathrm{SL}_2(\R)$
or $S_i = \mathrm{SL}_2(\R) \times \mathrm{SL}_2(\R)$, and
$V = V_1 \bigoplus V_2 \bigoplus ... \bigoplus V_m$
where each $V_i$ is an $L$-module such that each $S_j, j \neq i$, acts
trivially on $V_i$, $V_i$ is the standard representation of $\mathrm{SL}_2(\R)$
if $S_i = \mathrm{SL}_2(\R)$, and $V_i$ is the tensor product of two standard
representations of $\mathrm{SL}_2(\R)$ if $S_i = \mathrm{SL}_2(\R) \times \mathrm{SL}_2(\R)$.
\end{prop}

The next lemma is derived from \cite[Table 5, p.518]{[H]}.

\begin{lem}
\label{representation SL_2}
Let $Q$ be a simple real algebraic
group and $\mbox{rank}_{\R}Q \leq 1$. Let $d$ be the dimension of the minimal
nontrivial representation of $Q$ and $s$ be the dimension of the
symmetric space of $Q$. Then $d \geq s$ and $d = s$ if and only if $Q$
is isomorphic to $\mathrm{SL}_2(\R)$ and $d = s = 2$.
\end{lem}
{\bf Proof of Proposition \ref{representation}.} Let
$V = V_1 \bigoplus V_2 \bigoplus ... \bigoplus V_m$
be a direct sum of irreducible $L$-submodules. Each $V_i$ is a tensor
product of irreducible nontrivial $L_{ij}$-modules $V_{ij}$ (i.e. $V_i
= \bigotimes_{1 \le j \le r_i} V_{ij}, r_i \in {\N}$),
where $L_{ij} \in \{L_1$, $L_2,...,L_r\}$.
Put $n = \dim V$ and $n_{ij} = \dim
V_{ij}$. Then
\begin{equation}
\label{2}
n = \sum_{1\le i\le m} (\prod_{1\le j\le r_{i}}n_{ij}).
\end{equation}
For every $L_i$, we let $d_i$ (respectively, $s_i$) be the dimension of
the minimal nontrivial real representation of $L_i$ (respectively, the
dimension of the symmetric space of $L_i$). Remark that $s_1 + s_2 +
... + s_r$ is the dimension of the symmetric space of $L$. Using (\ref{1}), (\ref{2})
and the faithfulness of $\rho$,  we get
\begin{equation}
\label{3}
d_1 + d_2 + ... + d_r \le n \le s_1 + s_2 + ... + s_r.
\end{equation}
According to Lemma \ref{representation SL_2} $d_i \ge s_i$ for all $i$. It follows
from (\ref{3}) that $d_i = s_i = 2$ and $L_i = \mathrm{SL}_2(\R)$ for all $i$. In particular,
\begin{equation}
\label{4}
n = 2r.
\end{equation}
Since $\sum_{1\le i\le m} r_i \ge r$ and $n_{ij} \ge 2$ for all $i$
and $j$, it follows from (\ref{2}) and (\ref{4}) that  $\sum_{1\le i\le m} r_i =
r$, $1 \le r_i \le 2$ and all $n_{ij} = 2$ (i.e. each $V_{ij}$ is a standard
$\mathrm{SL}_2(\R)$-module). Moreover, we see that $V$ is a faithful representation of $L_1 \times \cdots \times L_r$.
This implies that $L = S_1 \times S_2 \times ... \times S_m$ as in the formulation of the proposition. \qed

\subsection{End of the proof.}
Let $\rho$ be as in the formulation of Proposition \ref{representation}. The
isomorphism $\mathcal{D}G/\mathrm{R}_u(\mathcal{D}G) \cong L$ gives a natural
surjective homomorphism $\pi: \mathcal{D}G \rightarrow L$. Put $\phi = \rho \circ \pi$.
Let $\gamma \in \mathcal{D}G$. By Proposition \ref{semisimple element} $\gamma_s$ fixes element-wise
a smooth curve on $G/H$. There exists $g \in \mathcal{D}G$ such that $g\gamma_sg^{-1} \in L$.
Hence $g\gamma_sg^{-1}$ fixes element-wise a smooth curve on $G/H$ passing through the origin.
So, $1$ is an eigenvalue of $\phi(g\gamma_sg^{-1})$ and, therefore, of $\phi(\gamma_s)$ and $\phi(\gamma)$ too. But $\phi(\Gamma)$
is Zariski dense in $L$. Therefore $1$ is an eigenvalue of $\rho(s)$ for every $s \in L$. In view of Proposition \ref{representation}
 $L$ is trivial, that is, $\Gamma$ is solvable. \qed

\section{On Auslander's conjecture}

\subsection{Some known results} 
First we formulate a general result which is often useful in tackling Auslander's conjecture.
So, let $S$ be a real, connected, non-compact, and semi-simple algebraic group. An element $g \in S$ is said to
be $\R$-regular if the number of eigenvalues having modulus $1$ (counted with multiplicity) of $\mathrm{Ad}(g)$
is minimal possible. Note that every $\R$-regular element in a semi-simple (or reductive) group is semi-simple.
It is known (and can be checked by direct computation)
that if $S = \mathrm{SL}_n(\R)$ or $\mathrm{Sp}_{2n}(\R), n \geq 2,$ (the cases arising in section 4.2) then
$g \in S$ is $\R$-regular if and only if all its eigenvalues are real
and their moduli are distinct. The following theorem is proved by different methods in \cite{[Be-L]}, \cite{[P]} and
\cite{[A-M-S1]}. (Concerning to its second part, we refer to \cite[Remark, p.545]{[P]}.)

\begin{thm}
\label{Prasad} Any Zariski dense sub-semigroup $\Delta$ of $S$ contains an $\R$-regular element. Moreover, the set of $\R$-regular elements
in $\Gamma$ is dense in $G$ in the Zariski topology.
\end{thm}

Note that if $S = \mathrm{SL}_n(\R)$ or $\mathrm{Sp}_{2n}(\R), n \geq 2,$ then the set of elements in $S$ with all eigenvalues different from $1$ is
Zariski open and non-empty which implies that the set of $\R$-regular elements
in $\Delta$ with all eigenvalues different from $1$ is Zariski dense in $G$. We will use this assertion in the course of
our proof in 4.2 of the Auslander conjecture for $n \leq 5$ \footnote{Note that Theorem \ref{Prasad}
is not indispensable for the proof of Auslander's conjecture for $\leq 5$. It could be replaced by a weaker claim
in the spirit of \cite[Lemmas 2.1-2.5]{[Ti]} but this would make the proof less natural and somewhat more complicate.}.

Further on, we denote by $\mathrm{SO}_{p,q}(\R)$ the special
orthogonal group of a quadratic form on $\R^n$ of signature $(p,q), \ n = p+q$, and by $\mathrm{Sp}_{2n}(\R)$ the
symplectic sub-group of $\mathrm{SL}_{2n}(\R)$. If $n = p$ we use the standard notation $\mathrm{SO}_{n}(\R)$ instead of
$\mathrm{SO}_{n,0}(\R)$.

Now, let $\lambda: \mbox{Aff}({\A}^n) \rightarrow \mbox{GL}_n({\R})$ be the natural
projection (see 1.2) and $H$ be  an
algebraic subgroup of $\GL_n(\R)$.
A subgroup $\Gamma \subset \mbox{Aff}({\A}^n)$
is called $H$-linear if $\lambda(\Gamma) \subset H$.  If $H = \mathrm{SO}_{n}(\R)$,
i.e. if $\Gamma$ consists of Euclidean transformations of ${\A}^n$, the
Auslander conjecture follows from the classical Bieberbach theorem.
Goldman and Kamishima (see [G-K]) proved the conjecture for
Lorentz space forms, i.e. for $H = \mathrm{SO}_{n-1,1}(\R)$, and
 Grunewald and Margulis proved it when
 $H$ is a reductive group of real rank $\leq 1$ (see [Gr-Mar]).

Recall the following results of Abels, Margulis and Soifer.

\begin{thm} \textsc{(}\cite[Theorems A and B]{[A-M-S3]}\textsc{)}
\label{AMS3} Suppose that $n = 2k + 1 \geq 3$. Then the following holds:
\begin{enumerate}
\item[(a)] if $k$ is even there is no $\Gamma$ acting properly discontinuously on $\A^n$ with
$\lambda(\Gamma)$ Zariski dense in $SO_{k+1,k}(\R)$;
\item[(b)] if $k$ is odd there are free groups $\Gamma$ acting properly discontinuously on $\A^n$ with
$\lambda(\Gamma)$ Zariski dense in $SO_{k+1,k}(\R)$.
\end{enumerate}
\end{thm}

Theorem \ref{AMS3}(a) was proved in the 1991 Margulis' preprint \cite{[Mar3]}. Theorem \ref{AMS3}(b)
is a generalization of Margulis' results \cite{[Mar1]} and \cite{[Mar2]} where Theorem \ref{AMS3}(b) is proved for $SO_{2,1}(\R)$
disproving a conjecture of J.Milnor \cite{[Mi]}. Different aspects of \cite{[Mar1]} and \cite{[Mar2]}
were developed in \cite{[Dr1]}, \cite{[Dr2]} and \cite{[Dr-G]}.

The results from \cite{[A-M-S3]} are sharpened by the following:

\begin{thm} \textsc{(}\cite[Theorems A,B and C]{[A-M-S4]}\textsc{)}
\label{AMS5} Suppose that $\lambda(\Gamma) \subset \mathrm{SO}_{p,q}(\R)$, $n = p+q$. Denote by $H$ the Zariski closure of $\lambda(\Gamma)$ in $\mathrm{GL}_n(\R)$. Then:
\begin{enumerate}
\item[(a)] $\Gamma$ can not act properly discontinuously on $\A^n$ if $|p-q| \geq 2$ and $H \supseteq \mathrm{SO}_{p,q}(\R)$;
\item[(b)] $\Gamma$ can not act properly discontinuously on $\A^n$ if $q$ is even and the homogeneous space $\mathrm{SO}_{p,q}(\R)/H$ is compact;
\item[(c)] $\Gamma$ is virtually solvable if $q = 2$ and $\Gamma \backslash {\A}^n$ is compact.
\end{enumerate}
\end{thm}

\subsection{Proof of Auslander's conjecture for $n \leq 5$} From now on $\Gamma$ is a subgroup of $\mathrm{Aff}(\A^n)$, $n \leq 5,$ and
$G$ is its Zariski closure in $\mathrm{Aff}(\A^n)$.
We will prove the following:

\begin{prop}
\label{small dimensions} If
$G$ contains a simple algebraic group of rank $\geq 2$ then $\Gamma$ does not act
properly discontinuously on $\A^n$.
\end{prop}

In view of Theorem \ref{Conj2} (or \cite{[To1]}), Proposition \ref{small dimensions} implies immediately:

\begin{thm}
\label{thm small dimensions}
The Auslander conjecture is true for $n \leq 5$.
\end{thm}

\subsubsection{}
We will use the notation and the terminology of section \ref{affine geometry}. In order to prove Proposition \ref{small dimensions} we need a particular case of the following general

\begin{prop}\textsc{(}cf.\cite[Lemma 4.2]{[To2]}\textsc{)}
\label{linear algebra}
 Let $H$ be a Zariski connected
algebraic subgroup of $\mathrm{Aff}({\A}^m)$, $S$ be a maximal reductive
subgroup of $H$ and $L$ be the Levi subgroup of $S$. Then there exists a decomposition
$V = W_1 \bigoplus \cdots \bigoplus W_k$,  where $W_i$ are irreducible $\lambda(S)$-modules, such that
for every $i \geq 1$ the sum $W_1 \bigoplus \cdots \bigoplus W_i$ is
$\lambda(H)$-invariant. Hence, for every $p \in \A^m$ there
 exists a frame with origin $p$ in which the matrix representation of $H$ is
of the form:
$$
\left(
\begin{array}{cccccccc}
{\rho_1} \\      \\        \\    \\  \\     \\ - \\ 0 \end{array}
\begin{array}{cccccccc}
      \\ {\rho_2}. \\      \\    \\ 0  \\    \\ - \\ 0 \end{array}
\begin{array}{cccccccc}
      \\          \\      \\ .  \\   \\      \\ - \\ . \end{array}
\begin{array}{cccccccc}
*      \\          \\      \\    \\ .  \\      \\ -  \\ . \end{array}
\begin{array}{cccccccc}
      \\          \\      \\    \\   \\ .      \\ -  \\ . \end{array}
\begin{array}{cccccccc}
      \\*      \\      \\    \\    \\ {\rho_k} \\ - \\ 0 \end{array}
\begin{array}{cccccccc}
 |    \\ |       \\   |   \\ | \\ |  \\ |      \\ | \\ | \end{array}
\begin{array}{cccccccc}
*     \\*      \\ .     \\ . \\ .  \\*    \\ - \\ 1 \end{array}
\right),
$$
where $\rho_i$ , $i=1,...,k$, are irreducible matrix representations of both $\lambda(H)$ and $\lambda(S)$.
Moreover, if some restriction ${\rho_i}_{\mid \lambda(L)}$ is not irreducible then
${\rho_i}_{\mid \lambda(L)} = \sigma_i \bigoplus \sigma_i $ where $\sigma_i$ is an
irreducible representation of $\lambda(L)$.
\end{prop}


{\bf Proof.}  The existence of the decomposition $V = W_1 \bigoplus \cdots \bigoplus W_k$ as in the formulation of the proposition will be proved by induction on $\dim V$.
Denote by $U$ the unipotent radical of $\lambda(H)$. Let $W_1$ be an irreducible $\lambda(H)$-submodule of $V$. Since $U$ is a unipotent group there exists
a $U$-invariant vector $\overrightarrow{v} \in W_1 \setminus \{0\}$. For every $g \in \lambda(S)$,  $g\overrightarrow{v}$
is also $U$-invariant. Since $\lambda(H) = \lambda(S)\ltimes U$, $W_1$ consists of $U$-invariant vectors and, therefore, is
an irreducible $\lambda(S)$-module. Suppose that $W_1, \cdots, W_i$ are irreducible nontrivial $\lambda(S)$-submodules of $V$ such that $W_1 + \cdots + W_i = W_1 \oplus \cdots \oplus W_i \subsetneqq V$ and $W_1 + \cdots + W_j$ is $\lambda(H)$-invariant for every $1 \leq j \leq i$. Let $W_{i+1}'$ be a $\lambda(H)$-submodule of $V$ containing $W_1 + \cdots + W_i$ and such that $W_{i+1}'/W_1 + \cdots + W_i$ is a non-trivial irreducible $\lambda(H)$-submodule of $V/W_1 + \cdots + W_i$. By the complete reducibility of the action of $\lambda(S)$ on $W_{i+1}'$ there exists an irreducible $\lambda(S)$-submodule $W_{i+1}$ of $W_{i+1}'$ such that $W_{i+1}' = W_1 \oplus \cdots \oplus W_{i+1}$, completing the proof of the existence of the decomposition.

Now, suppose that $\rho: S \rightarrow \GL(W)$ is an irreducible representation but $\rho_{\mid L}$ is not. Let $Z$ be the center of $S$.
If $W'$ is an irreducible $Z$-module then $\dim W' = 2$ and $W$ is a direct sum of translations of $W'$ by elements from $S$. This
implies that if $W''$ is an irreducible $L$-module then there exists a $c \in Z$ such that $W = W'' \oplus cW''$, proving the last
assertion of the proposition. \qed

\medskip

\textit{Remark.} It is worth mentioning that $\lambda(S)$ is isomorphic to $S$ and the representations $\{{\rho_i}_{\mid \lambda(S)}: i = 1, \cdots, k\},$ do not depend (up to isomorphism)
on the choice of $S$.
\subsubsection{Proof of Proposition \ref{small dimensions}.} We divide the proof in two steps: the first being mostly routine and the second
containing the main ingredients of the proof.

\textit{Step 1.}
Replacing $\Gamma$ by a subgroup of finite index we may (and will) assume that
the algebraic group $G$ is Zariski connected. We will apply Proposition \ref{linear algebra} with $H = G$.
 We fix a frame $\{p; \overrightarrow{v_1}, \cdots, \overrightarrow{v_n}\}$ of
$\A^n, n \leq 5,$ such that the matrix representation of $G$ in this frame is as given by Proposition \ref{linear algebra} and we keep the notation $S$, $L$
and $\rho_i, \ i = 1, \cdots, k,$ from its formulation. Since $\lambda_{\mid S}$ is injective, when this does not lead to confusion, we also write $\rho_i$ instead of $\rho_i \circ \lambda$.

There exists $l$ such that $\mathcal{D}^{l}G = \mathcal{D}^{l+1}G$. Since $\mathcal{D}^{l}\Gamma$
is Zariski dense in $\mathcal{D}^{l}G$ and $\mathcal{D}^{l}G$ contains a simple algebraic group of rank $\geq 2$ (because $G$ does),
 replacing $\Gamma$ by $\mathcal{D}^{l}\Gamma$, we reduce the proof of Proposition \ref{small dimensions} to the case when $G = \mathcal{D}G$. In this case $L = S$. In view of the classification in \cite[Table 2]{[Bou]} of the dimensions of the irreducible representations of the simple algebraic groups, if $H$ is a simple algebraic group of rank $d > 0$ and $\rho: H \rightarrow \mathrm{SL}(W)$ is its non-trivial representation then
$\dim W \geq d+1$ and $\dim W = d+1$ if and only if $\rho$ is an isomorphism \footnote{Formally, \cite[Table 2]{[Bou]} concerns the representations
of the algebraic groups over $\C$ but since the absolute rank (over $\C$) of an algebraic group over $\R$ is greater than or equal to its real rank the assertion remains true for real
algebraic groups as in the context of this paper.}. It follows from Proposition \ref{linear algebra} and the assumption that $L$ contains a simple group of rank $\geq 2$ that: (a) if $n = 3$ then $L = \mathrm{SL}_{3}(\R)$, and (b) if $L$ is not simple then $n = 5$ and $L = \mathrm{SL}_{3}(\R) \times \mathrm{SL}_{2}(\R)$.
Using Theorem \ref{Prasad}, we get that in both cases (a) and (b) there exists an element $\gamma \in \Gamma$  of infinite order such that all eigenvalues of $\lambda(\gamma)$ are different from 1 and pairwise different. According to Proposition \ref{semisimple element}(iii), $\gamma$
admits a fixed point proving that $\Gamma$ does not act properly discontinuously on $\A^n$.

So, it remains to consider the case when $4 \leq n \leq 5$ and $L$ is a simple algebraic group of rank $\geq 2$.
It is enough to prove that if Proposition \ref{small dimensions} is valid for $n-1$ it is also valid for $n$. Assume to the contrary that $\Gamma$ acts properly discontinuously on $\A^n$ and Proposition \ref{small dimensions} is true for $n-1$. If one of the representations $\rho_i, 1 \leq i \leq k,$ as in the formulation of Proposition \ref{linear algebra}, is not trivial then $\dim \rho_i \geq 3$
and, since $n \leq 5$, all remaining representations $\rho_j, j \neq i$, are trivial. Also, if $\rho_k$ is trivial it follows from
 $G = \mathcal{D}G$ that $\Gamma$ acts properly discontinuously of the hyperplane $p + \R\overrightarrow{v_1} + \cdots + \R\overrightarrow{v}_{n-1}$ contradicting the assumption that Proposition \ref{small dimensions} is true for $n-1$. Hence $\rho_k$ is a non-trivial representation, $3 \leq \dim \rho_k \leq 5$, and all $\rho_i$, $1 \leq i \leq k-1$, are trivial. Moreover, it follows from the remark after Proposition \ref{semisimple element}
 that $\lambda(L)$ should not contain an element with pairwise different eigenvalues all different from $1$. Now, let $k = 1$. In view of \cite[Table 2]{[Bou]}, $L = \mathrm{SO}_{3,2}(\R)$, $n = 5$ and ${\rho_1}_{|L}$ is the standard representation.
 By Margulis' Theorem \ref{AMS3}(a), in this case the action of $\Gamma$ on $\A^5$ is not properly discontinuous\footnote{The treatment of this case in  \cite{[To2]} contains an error.}.
Therefore $k > 1$ and $3 \leq \dim \rho_k \leq 4$. Since $\rank (L) \geq 2$, \cite[Table 2]{[Bou]} implies that $\rho_k$ is one of the standard
representations of $\mathrm{SL}_3(\R)$, $\mathrm{SL}_4(\R)$ or $\mathrm{Sp}_4(\R)$. So, in order to complete the proof it remains to consider the
possibilities $\mathrm{(i)}$, $\mathrm{(ii)}$, $\mathrm{(iii)}$ and $\mathrm{(iv)}$ below.
(In the formulations of $\mathrm{(i)} - \mathrm{(iv)}$ we use the notation: $\mathrm{M}_{i,j}(\R)$ is the set of all real matrices with $i$ rows and
$j$ columns, $0_{i,j}$ is the zero in $\mathrm{M}_{i,j}(\R)$ and $I_i$ is the unit matrix in $\mathrm{M}_{i,i}(\R)$.)

\medskip

  $\mathrm{(i)}$ $\lambda(G) \subseteq \Big\{ \left(\begin{array}{cc} \mathrm{I}_1&m_{1,4}\\ 0_{4,1}&g\\ \end{array}
\right)\mid g \in L, m_{1,4} \in \mathrm{M}_{1,4}(\R) \Big\}$ and $L = \mathrm{SL}_4(\R)$;

$\mathrm{(ii)}$ $\lambda(G) \subseteq \Big\{ \left(\begin{array}{cc} \mathrm{I}_1&m_{1,4}\\ 0_{4,1}&g\\ \end{array}
\right)\mid g \in L, m_{1,4} \in \mathrm{M}_{1,4}(\R) \Big\}$ and $L = \mathrm{Sp}_4(\R)$;

 $\mathrm{(iii)}$ $\lambda(G) \subseteq \Big\{\left(\begin{array}{cc}\mathrm{I}_2&m_{2,3}\\ 0_{3,2}&g\\ \end{array}
\right)\mid g \in L, m_{2,3} \in \mathrm{M}_{2,3}(\R) \Big\}$ and $L=\mathrm{SL}_3(\R)$;

$\mathrm{(iv)}$ $\lambda(G) \subseteq \Big\{ \left(\begin{array}{cc} \mathrm{I}_1&m_{1,3}\\ 0_{3,1}&g\\ \end{array}
\right)\mid g \in L, m_{1,3} \in \mathrm{M}_{1,3}(\R) \Big\}$ and $L = \mathrm{SL}_3(\R)$.

\medskip

\textit{Step 2.}
We treat the cases $\mathrm{(i)} - \mathrm{(iv)}$ simultaneously.
By Theorem \ref{Prasad} there exists $\gamma \in \Gamma$ such that $\rho_k(\gamma)$ is $\R$-regular in
the Zariski closure of $\rho_k(\Gamma)$ and all its eigenvalues are different from $1$. (Note that this property
of $\gamma$ remains valid for any choice of $L$ in \textit{Step 1}.)
Let $\gamma = \gamma_s \gamma_u$ be the Jordan decomposition of $\gamma$. Choose $L$ such that
$\gamma_s \in L$. 
Fix a $p_\circ \in \A^n$ with $Lp_\circ = p_\circ$.

We have $V = W_\circ \oplus W$ where $W_\circ = \{\overrightarrow{x} \in V | \lambda(\gamma_s)(\overrightarrow{x}) = \overrightarrow{x}\} = \{\overrightarrow{x} \in V | \lambda(g)(\overrightarrow{x}) = \overrightarrow{x} \ \textmd{for all} \ g \in G\}$ and $W$ is $\lambda(\gamma_s)$-invariant. Using $\lambda(\gamma_s) \lambda(\gamma_u) = \lambda(\gamma_u)\lambda(\gamma_s)$ and the choice of $\gamma$, one proves easily that $\gamma_u$ is a translation belonging to the center of $G$.
Indeed, let $\lambda(\gamma) = \left(\begin{array}{cc} \mathrm{I}_i&m_{i,j}\\ 0_{j,i}&g\\ \end{array} \right)$ and $g = g_s g_u$ be the Jordan decomposition of $g$. By the choice of $L$ and Proposition \ref{linear algebra}, we have $\lambda(\gamma_s) = \left(\begin{array}{cc} \mathrm{I}_i&0_{i,j}\\ 0_{j,i}&g_s\\ \end{array} \right)$ and $\lambda(\gamma_u) = \left(\begin{array}{cc} \mathrm{I}_i&m_{i,j}\\ 0_{j,i}&g_u\\ \end{array} \right)$.
Since the eigenvalues of $g_s$ are pairwise different we get that $g_u = I_j$ and since they are all different from $1$ we get that $m_{i,j} = 0_{i,j}$. Therefore, $\lambda(\gamma_u) = \mathrm{Id}_V$, equivalently, $\gamma_u$ is a translation by a vector $\overrightarrow{v_\gamma}$. We have
$$
\gamma_s \gamma_u p_\circ = p_\circ + \lambda(\gamma_s)(\overrightarrow{v_\gamma}) = \gamma_u \gamma_s p_\circ = p_\circ + \overrightarrow{v_\gamma}.
$$
Hence $\overrightarrow{v_\gamma} \in W_\circ$. Let $h \in G$. Then $\lambda(h)\overrightarrow{v_\gamma} = \overrightarrow{v_\gamma}$. So, if $x  \in \A^n$ then
$$
h \gamma_u(x) = h(x + \overrightarrow{v_\gamma}) = h(x) + \overrightarrow{v_\gamma} = \gamma_u h(x),
$$
proving that $\gamma_u$ belongs to the center of $G$.

Put
$$
E^{\circ}(\gamma) = \{p \in \A^n \mid \gamma_s p = p \}.
$$
Then $E^{\circ}(\gamma) = p_\circ + W_\circ$.
Denote
$$
V^+(\gamma) = \{\overrightarrow{v} \in V \mid \underset{n \rightarrow -\infty}{\lim}\gamma_s^n(p_\circ + \overrightarrow{v})  =  p_\circ\}
$$
and
$$
V^-(\gamma) = \{\overrightarrow{v} \in V \mid \underset{n \rightarrow +\infty}{\lim}\gamma_s^n(p_\circ + \overrightarrow{v})  =  p_\circ\}.
$$
Let
$$
E^{+}(\gamma) = E^{\circ}(\gamma) + V^+(\gamma) \ \ \text{and} \ \ E^{-}(\gamma) = E^{\circ}(\gamma) + V^-(\gamma).
$$
Since $E^{\pm}(\gamma) = E^{\mp}(\gamma^{-1})$, we may (and will) assume that
$$\dim E^{+}(\gamma) \geq \dim E^{-}(\gamma).$$

Let $\delta \in \Gamma$. Since $\gamma_u$ is central, we have $\gamma_u = (\delta\gamma \delta^{-1})_u$. Also,
$$
E^{\circ}(\delta\gamma \delta^{-1}) = \delta E^{\circ}(\gamma)  \ \ \text{and} \ \ E^{\pm}(\delta\gamma \delta^{-1}) = \delta E^{\pm}(\gamma).
$$
Remark that $E^{\circ}(\gamma)$ and $E^{\circ}(\delta\gamma \delta^{-1})$ are parallel subspaces of $\A^n$ directed by $W_\circ$. Suppose that
$E^{\circ}(\gamma) = \delta E^{\circ}(\gamma)$ for all $\delta \in \Gamma$.
 Recall that every subgroup of $\mathrm{Aff}(\A^m), m \leq 2,$ acting properly discontinuously on $\A^m$ is virtually solvable. Since $G$ is connected, $G = \mathcal{D}G$ and $\dim E^{\circ}(\gamma) \leq 2$, we get that
the action of $\Gamma$ on $E^{\circ}(\gamma)$ is trivial which contradicts the assumption that
$\Gamma$ acts properly discontinuously on $\A^n$.
Therefore there exists $\delta \in \Gamma$ such that
\begin{equation}
\label{5}
E^{\circ}(\gamma) \cap E^{\circ}(\delta\gamma \delta^{-1}) = \emptyset.
\end{equation}
With such a $\delta$, $E^{+}(\gamma) \cap E^{+}(\delta\gamma \delta^{-1})$ contains a line $l = q + \R \overrightarrow{v_\gamma}$.
We can write $q = q_1 + \overrightarrow{w_1}$, where $q_1 \in E^{\circ}(\gamma)$ and $\overrightarrow{w_1} \in V^+(\gamma) \setminus \{\overrightarrow{0}\}$,
and $q = q_2 + \overrightarrow{w_2}$, where $q_2 \in E^{\circ}(\delta\gamma \delta^{-1})$ and $\overrightarrow{w_2} \in V^+(\delta\gamma \delta^{-1}) \setminus \{\overrightarrow{0}\}$.

Put $p_m = q + m\overrightarrow{v_\gamma}, m \in \N$.
Then
$$
\underset{m \rightarrow +\infty}{\lim}\gamma^{-m}(p_m)  = q_1  \ \ \text{and} \ \ \underset{m \rightarrow +\infty}{\lim}(\delta\gamma \delta^{-1})^{-m}(p_m)  = q_2.
$$
Note that the set $X = \{\gamma^{-m}(p_m) \mid m \in \N \} \cup \{(\delta\gamma \delta^{-1})^{-m}(p_m) \mid m \in \N \}$ is relatively compact
and $(\delta\gamma^{-m} \delta^{-1}\gamma^{m})X \cap X \neq \emptyset$ for all $m$. It follows from (\ref{5}) that  $\delta\gamma^{-m_1} \delta^{-1}\gamma^{m_1} \neq \delta\gamma^{-m_2} \delta^{-1}\gamma^{m_2}$ if $m_1 \neq m_2$ which contradicts the assumption that $\Gamma$ acts
properly discontinuously on $\A^n$. This completes our proof. \qed

\vskip.5cm

{\it Remark:} As noted in \cite{[To2]}, the element $\delta$ can be chosen in such a way that the subgroup
spanned by $\gamma^m$ and $\delta\gamma^{m}\delta^{-1}$ is free if $m$ is sufficiently large. This can be achieved by a well-known argument of Tits \cite{[Ti]}. Indeed, let $\varphi: G \rightarrow L$ be the natural projection and $L$ be identified with its image in $\mathrm{SL}(W)$. Let $\mathrm{P}(W)$ be the projective space of $W$. If $g \in G$ and $\varphi(g)$ has pairwise different
positive eigenvalues $\alpha_1 > \cdots > \alpha_t > 0$ with respective eigenvectors $\overrightarrow{w_1}, \cdots, \overrightarrow{w_t}$,
we denote by $\mathrm{A}(g)$ (resp.$\mathrm{A}'(g))$ the projectivization of the vector space $\R\overrightarrow{w_1}$ (resp. $\R\overrightarrow{w_2} + \cdots + \R\overrightarrow{w_t}$). Since $L$ acts irreducibly on $W$, we can choose
$\delta$ in such a way that $\mathrm{A}(\gamma) \cup \mathrm{A}(\gamma^{-1}) \subset \mathrm{P}(W) \setminus (\mathrm{A}'(\delta\gamma \delta^{-1}) \cup \mathrm{A}'(\delta\gamma^{-1}\delta^{-1}))$ and $\mathrm{A}(\delta\gamma \delta^{-1}) \cup \mathrm{A}(\delta\gamma^{-1}\delta^{-1}) \subset \mathrm{P}(W) \setminus (\mathrm{A}'(\gamma) \cup \mathrm{A}'(\gamma^{-1}))$. Now, it follows from the "ping-pong lemma" that the
subgroup spanned by $\gamma^m$ and $\delta\gamma^{m}\delta^{-1}$ is free if $m$ is sufficiently large.

\medskip

\textit{Acknowledgements:}  After the initial version of this paper has been written several people contributed to its improvement.
Thanks are due to Grisha Margulis and Andrei Rapinchuk for informing us about the existence of the papers \cite{[A-M-S5]} and \cite{[D-P]}, respectively, and to Gopal Prasad for his valuable remarks and suggestions.

\end{document}